\documentclass[reqno,centertags]{amsart}
\usepackage{amsmath,amsthm,amscd,amssymb}
\usepackage{latexsym}
\newcommand{\bbC}{{\mathbb{C}}}
\newcommand{\bbD}{{\mathbb{D}}}

\newcommand{\bbG}{{\mathbb{G}}}

\newcommand{\bbR}{{\mathbb{R}}}

\newcommand{\lb}{\label}
\newcommand{\f}{\frac}

\newcommand{\ol}{\overline}

\newcommand{\wti}{\widetilde  }

\newcommand{\tr}{\text{\rm{Tr}}}

\newcommand{\Sz}{\text{\rm{Sz}}}

\newcommand{\ran}{\text{\rm{ran}}}

\newcommand{\s}{\text{\rm{s}}}

\newcommand{\bi}{\bibitem}

\newcommand{\beq}{\begin{equation}}
\newcommand{\eeq}{\end{equation}}
\newcommand{\ba}{\begin{align}}
\newcommand{\ea}{\end{align}}
\newcommand{\veps}{\varepsilon}

\newcounter{smalllist}
\newenvironment{SL}{\begin{list}{{\rm\roman{smalllist})}}{%
\setlength{\topsep}{0mm}\setlength{\parsep}{0mm}\setlength{\itemsep}{0mm}%
\setlength{\labelwidth}{2em}\setlength{\leftmargin}{2em}\usecounter{smalllist}%
}}{\end{list}}

\DeclareMathOperator{\Ima}{Im}

\allowdisplaybreaks
\numberwithin{equation}{section}

\newtheorem{theorem}{Theorem}[section]

\newtheorem*{p2.1}{Proposition 2.1}
\newtheorem{proposition}[theorem]{Proposition}
\newtheorem{lemma}[theorem]{Lemma}

\theoremstyle{definition}
\newtheorem{example}[theorem]{Example}

\theoremstyle{remark}
\newtheorem*{remark}{Remark} 
\newtheorem*{remarks}{Remarks}

\theoremstyle{definition}
\newtheorem*{definition}{Definition}

\newcommand{\abs}[1]{\lvert#1\rvert}

\begin{document}
\title[Meromorphic Jost Functions]{Jost Functions and Jost Solutions for Jacobi Matrices,
III.~Asymptotic Series for Decay and Meromorphicity}
\author[B.\ Simon]{Barry Simon$^*$}

\thanks{$^*$ Mathematics 253-37, California Institute of Technology, Pasadena, CA 91125.
E-mail: bsimon@caltech.edu. Supported in part by NSF grant DMS-0140592 and in part 
by Grant No.\ 2002068 from the United States-Israel Binational Science Foundation 
(BSF), Jerusalem, Israel}

\date{March 7, 2005} 

\begin{abstract}  We show that the parameters $a_n,b_n$ of a Jacobi matrix have a 
complete asymptotic series 
\begin{align*}
a_n^2 -1 &= \sum_{k=1}^{K(R)} p_k(n) \mu_k^{-2n} + O(R^{-2n}) \\ 
b_n &= \sum_{k=1}^{K(R)} p_k(n)\mu_k^{-2n+1} + O(R^{-2n}) 
\end{align*} 
where $1<\abs{\mu_j}<R$ for $j\leq K(R)$ and all $R$ if and only if the Jost function, 
$u$, written in terms of $z$ (where $E=z+z^{-1}$) is an entire meromorphic function. 
We relate the poles of $u$ to the $\mu_j$'s.
\end{abstract}

\maketitle

\section{Introduction} \lb{s1} 

In this paper, we are going to consider semi-infinite Jacobi matrices 
\begin{equation} \lb{1.1} 
J=\begin{pmatrix} b_1 & a_1 & 0 & \dots \\
a_1 & b_2 & a_2 & \dots \\
0 & a_2 & b_3 & \dots \\
\vdots & \vdots & \vdots &  \ddots
\end{pmatrix} 
\end{equation}  
whose Jacobi parameters have exponential decay (i.e., $\limsup_{n\to\infty} (\abs{b_n} 
+ \abs{a_n-1})^{1/n} <1$. As explained in the first two papers of this series 
\cite{Jost1,Jost2} (and well-known earlier), such a $J$ has an associated Jost 
function, $u$, defined and analytic in a neighborhood of $\ol{\bbD}$ where $\bbD= 
\{z\in\bbC \mid\abs{z}<1\}$. 

As is standard, $J$ describes the recursion relations for orthogonal polynomials on the 
real line (OPRL). There is a probability measure, $\gamma$, so that the orthonormal 
polynomials, $p_n(x)$ \cite{Szb,OPUC1,OPUC2}, defined by $\gamma$ obey 
\begin{equation} \lb{1.2} 
xp_n(x) = a_{n+1} p_{n+1}(x) + b_{n+1} p_n(x) + a_n p_{n-1}(x) 
\end{equation} 
$\gamma$ is the spectral measure for $J$ and vector $(1\,0\, 0 \dots)^t$ and the $a$'s 
and $b$'s can be obtained from $\gamma$ by Gram-Schmidt on the moments. 

$u$ is defined by $\gamma$ via the following three facts: 
\begin{SL} 
\item[(i)] $u(z)=0$ for $z\in\bbD$ if and only if $z+z^{-1}$ is an eigenvalue of $J$\!. 
\item[(ii)] The support of $d\gamma_\s$, the singular part of $d\gamma$, is a finite set 
of eigenvalues in $\bbR\backslash [-2,2]$ and 
\begin{equation} \lb{1.3} 
d\gamma\restriction [-2,2] = f(x)\, dx 
\end{equation} 
where for any $\theta\in [0,2\pi)$,  
\begin{equation} \lb{1.4} 
f(2\cos\theta) = \f{1}{\pi} \biggl[ \f{\sin\theta}{\abs{u(e^{i\theta})}^2}\biggr]  
\end{equation} 
\item[(iii)] $u(0)>0$. 
\end{SL} 

These determine $u$ by standard theory of nice analytic functions on $\bbD$. $u$ does 
not determine $\gamma$ in many cases. For by \eqref{1.3}/\eqref{1.4}, $u$ determines 
the a.c.\ part of $\gamma$ and the positions of the pure points but not their weights. 
We prefer to normalize the weights by looking at 
\begin{equation} \lb{1.5} 
M(z) = -\int \f{d\rho(x)}{x-(z+z^{-1})}  
\end{equation} 
and looking at the residues of $M$ at the points where $u(z)=0$. 

In any event, $u$ plus the weights are spectral data, and our goal here is to produce 
equivalences between this spectral data side and the recursion coefficient side. 

To state our main theorems, we define 

\begin{definition} A sequence, $(x_0, \dots, x_n, \dots)$, of complex numbers is said to have 
an asymptotic series up to $R>1$ if and only if there exists $\mu_1, \dots, \mu_{K(R)}$ 
in $\{z\mid 1<\abs{z}<R\}$ and polynomials $p_1, \dots, p_{K(R)}$ so that 
\begin{equation} \lb{1.6} 
\limsup_{n\to\infty}\, \bigg| x_n -\sum_{j=1}^{K(R)} p_j(n) \mu_j^{-n}\biggr|^{1/n} 
\leq R^{-1} 
\end{equation} 
We say $(x_0, \dots)$ has a complete asymptotic series if it has one for each $R>1$. 
\end{definition} 

It is easy to see that the $x$'s uniquely determine the $p$'s and $\mu$'s and that 

\begin{theorem}\lb{T1.1} $\{x_n\}_{n=0}^\infty$ has an asymptotic series up to $R$ if 
and only if  
\begin{equation} \lb{1.7} 
f(z)\equiv \sum_{n=0}^\infty x_n z^n  
\end{equation} 
is meromorphic in $\{z\mid \abs{z} <R\}$ with no singularities in a neighborhood of 
$\ol{\bbD}$ and finitely many poles in the region. $\{x_n\}_{n=1}^\infty$ has a complete 
asymptotic series if and only if $f$ is entire meromorphic. 
\end{theorem} 

Indeed, the poles are at the $\mu_j$ and their orders are one plus the degrees of the $p_j$. 

We say a set of Jacobi parameters has an asymptotic series up to $R$ if and only if the 
sequence 
\begin{equation} \lb{1.8} 
(1,-b_1, 1-a_1^2, -b_2, 1-a_2^2, \dots)
\end{equation} 
has an asymptotic series up to $R$. Thus, the function $f$ is 
\begin{equation} \lb{1.9} 
B(z)=1-\sum_{n=0}^\infty\, [b_{n+1} z^{2n+1} + (a_{n+1}^2 -1) z^{2n+2}] 
\end{equation} 

$B$ will enter naturally below, but we note the following interpretation: If $J_0$ is the 
Jacobi matrix with $a_n=1$, $b_n\equiv 0$, and $\delta J=J-J_0$, then (see Lemma~6.2 of 
\cite{Jost1}): 
\begin{equation} \lb{1.10} 
\tr(\delta J(J_0-(z+z^{-1}))^{-1}) = -(z^{-1} -z)^{-1} 
\biggl\{\,\sum_{n=1}^\infty b_n (1-z^{2n}) + 2\sum_{n=1}^\infty (a_n -1) 
(z-z^{2n+1})\biggr\}
\end{equation} 
Moreover (see Theorem~2.16 of \cite{KS}), 
\begin{equation} \lb{1.11} 
u(z)=\biggl(\,\prod_{j=1}^\infty a_j\biggr)^{-1} \det (1+\delta J(J_0 -(z+z^{-1}))^{-1}) 
\end{equation}  

Taking into account that $a_n^2 -1=2(a_n -1) + O((a_n-1)^2)$ and $\det (1+A)=1 + 
\tr(A) + O(\|A\|_1^2)$, we see that if $\delta J$ is trace class, then 
\begin{equation} \lb{1.12} 
-(z^{-1}-z) \biggl(\,\prod_{j=1}^\infty a_j\biggr) u(z) = c+ zB(z) + O(\|\delta J\|_1^2)  
\end{equation} 
for a constant, $c$. Thus, $B(z)$ is a kind of first-order (Born) approximation to $u$.  

In some ways, our main result in this paper is 

\begin{theorem}\lb{T1.2} The Jacobi parameters have a complete asymptotic series if and only 
if $u$ is an entire meromorphic function. Equivalently, $B(z)$ is entire meromorphic if and 
only if $u(z)$ is.
\end{theorem}  

Of course, one wants to understand the relation between the poles of $u$ and those of $B$. 
Both for that understanding and because we will actually use them in our proofs in 
Section~\ref{s3}, it pays to review our recent results \cite{MSF} on the analogous 
problem for orthogonal polynomials on the unit circle (OPUC). The basics (see 
\cite{OPUC1,OPUC2} for background) associate to a nontrivial probability measure, 
$\mu$, on $\partial\bbD$ a sequence of Verblunsky coefficients defined by 
\begin{equation} \lb{1.13} 
\Phi_{n+1}(z) = z\Phi_n(z) - \bar\alpha_n \Phi_n^*(z) 
\end{equation} 
where $\Phi_n$ are the monic orthogonal polynomials for $\mu$ and 
\begin{equation} \lb{1.14} 
\Phi_n^*(z) = z^n \, \ol{\Phi_n (1/\bar z)} 
\end{equation} 
In place of $B$, \cite{MSF} uses 
\begin{equation} \lb{1.15} 
S(z)=1-\sum_{j=1}^\infty \alpha_{j-1} z^j 
\end{equation} 
and, in place of $u$, the Szeg\H{o} function 
\begin{equation} \lb{1.16} 
D(z)=\exp \biggl( \int \f{e^{i\theta}+z}{e^{i\theta}-z} \, \log(w(\theta))\, 
\f{\,d\theta}{4\pi\,}\biggr) 
\end{equation}  
where $d\mu = w(\theta)\f{d\theta}{2\pi} + d\mu_\s$. One also defines 
\begin{equation} \lb{1.17} 
r(z)=\f{D^{-1}(z)}{\,\ol{D^{-1} (1/\bar z)}\,} 
\end{equation} 
The main theorems of \cite{MSF} are: 

\begin{theorem}[\cite{DO}]\lb{T1.3} If $\limsup\abs{\alpha_n}^{1/n}=R^{-1} <1$, 
then $r(z)-S(z)$ is analytic in $\{z\mid 1-\delta < \abs{z} < R^3\}$ for some 
$\delta >0$.  
\end{theorem} 
 
\begin{remarks} 1. This result is due to Deift-Ostensson \cite{DO}, but \cite{MSF} has 
a new proof. Earlier, \cite{OPUC1} proved the weaker result when $R^3$ is replaced 
by $R^2$. 

2. The point is that $r$ and $S$ both have singularities on $\abs{z}=R$. This theorem 
says they cancel, as do other singularities in $\{z\mid R<\abs{z} < R^3\}$. 

3. \cite{MSF} has explicit examples where $r(z)-S(z)$ has singularities on $\{z\mid 
\abs{z}=R^3\}$ and shows that this is the case generically. So $R^3$ is best possible. 
\end{remarks}

Given a discrete set, $\Omega\subset \{z\mid\abs{z}>1\}$, with limit points only at 
$\infty$, we define 
\begin{align} 
\bbG^{2j-1}(\Omega) &= \{\mu_1 \dots \mu_j\bar\mu_{j+1}\dots\bar\mu_{2j-1}\mid 
\mu_k\in\Omega\} \lb{1.18} \\
\bbG(\Omega) &=\bigcup_{j=1}^\infty \bbG^{2j-1}(\Omega)  \lb{1.19}
\end{align} 

\begin{theorem}[\cite{MSF}] \lb{T1.4} $S$ is entire analytic if and only if $D^{-1}$ is. 
If $T$ is the set of poles of $S(z)$ and $P$ the poles of $D^{-1}(z)$, then 
\begin{equation} \lb{1.20} 
T\subset \bbG(P) \qquad P\subset \bbG(T) 
\end{equation} 
\end{theorem} 

Analogously to Theorem~\ref{T1.3}, we will prove in Section~\ref{s2} that 

\begin{theorem}\lb{T1.5} Suppose 
\begin{equation} \lb{1.21} 
\limsup_{n\to\infty}\, (\abs{a_n^2 -1} + \abs{b_n})^{1/2n} = R^{-1} <1 
\end{equation} 
Then 
\begin{equation} \lb{1.22} 
(1-z^2) u(z) +z^2 \, \ol{u(1/\bar z)}\, B(z) 
\end{equation}  
is analytic in $\{z\mid R^{-1} < \abs{z} <R^2\}$. 
\end{theorem} 

\begin{remarks} 1. \cite{Jost2} has necessary and sufficient conditions on $\{u, 
\text{weights}\}$ for \eqref{1.21} to hold. If there are no eigenvalues of $J$ outside 
$[-2,2]$, the condition is that $u$ is analytic in $\{z\mid\abs{z} <R\}$. 

2. $u$ is real on $\bbR$ so $\ol{u(\bar z)} =u(z)$ and thus, \eqref{1.22} could be 
written $(1-z^2) u(z) +u(1/z) B(z)$; we write it as we do for analogy with the OPUC case. 

3. The point, of course, is that $B$ has singularities on $\{z\mid\abs{z}=R\}$, so 
this theorem implies a cancellation either via zeros of $\ol{u(1/\bar z)}$ or singularities 
of $u$. Since $\ol{u(1/\bar z)}$ can have zeros in $\abs{z}>1$ (while $\ol{D(\bar z)^{-1}}$ 
cannot), the situation is somewhat different from OPUC. We will discuss this further in 
Section~\ref{s2}. 

4. As we will show in Section~\ref{s3}, the function in \eqref{1.22} often has a singularity 
at $z=R^2$, so one cannot increase the $R^2$ to $R^3$ as one can in the OPUC case. 
The reason for this difference will become clear in Section~\ref{s3}. 
\end{remarks} 

For the analog of Theorem~\ref{T1.4}, we need to define a larger set than $\bbG$. 
In our situation, $u$ and $B$ are real on $\bbR$ so their poles are symmetric about 
$\bbR$. So for this, we will suppose $\Omega\subset\{z\mid\abs{z}>1\}$ with limit 
point only at infinity, and
\begin{equation} \lb{1.23} 
\!\ol{\Omega} =\Omega 
\end{equation} 
In that case, for any $m$, we define 
\begin{equation} \lb{1.24} 
\bbG^{(m)} (\Omega) =\{\mu_1 \dots \mu_m\mid \mu_k\in\Omega\} 
\end{equation} 
When \eqref{1.23} holds, this agrees with the previous definition if $m=2k-1$, 
\begin{equation} \lb{1.25} 
\wti\bbG(\Omega) = \biggl[\, \bigcup_{m=1}^\infty \bbG^{(m)} (\Omega) \biggr] 
\cup \biggl[ -\bigcup_{m=1}^\infty \bbG^{(m)} (\Omega)\biggr]
\end{equation}  

Our main results refine Theorem~\ref{T1.2}: 

\begin{theorem}\lb{T1.6} Let $J$ have no spectrum outside $[-2,2]$ and let $u$ be entire 
meromorphic and nonvanishing at $z=\pm 1$. Let $P$ be the poles of $u$ and $T$ the 
poles of $B$. Then 
\begin{equation} \lb{1.26} 
P\subset\wti\bbG(T) \qquad T\subset\wti\bbG(P) 
\end{equation}  
\end{theorem} 

To state the result when there are bound states, we recall and extend a notion from 
\cite{Jost2}. 

\begin{definition} Let $u$ be a meromorphic function and $z_0\in\bbD$ a point with 
$u(z_0)=0$ (so $z_0$ is real and $z_0 + z_0^{-1}\in\sigma (J)$). $z_0$ is called 
a noncanonical zero for $J$ if and only if $1/z_0$ is not a pole of $u$ and 
\begin{equation} \lb{1.27}  
\lim_{z\to z_0}\, (z-z_0) M(z) \neq -(z_0-z_0^{-1}) 
\biggl[ u'(z_0) u\biggl(\f{1}{z_0}\biggr)\biggr]^{-1}
\end{equation} 
\end{definition} 

Thus, $z_0$ is not noncanonical (which we will call canonical) if $u$ is regular 
at $1/z_0$ and equality holds in \eqref{1.27}. Here is what we will prove in case 
there are bound states or $u(\pm 1)=0$: 

\begin{theorem}\lb{T1.7} Suppose $u$ is entire meromorphic. Let $T$ be the poles 
of $B$. Let $P_1$ be the poles of $u$ and $P_2$ the $\{z^{-1}\mid z$ is a 
noncanonical zero for $J\}$. Let $P=P_1\cup P_2$. Then \eqref{1.26} holds. 
\end{theorem} 

As in \cite{MSF}, one can easily prove results relating meromorphicity of $u$ in 
$\{z\mid\abs{z} <R^{2\ell-1}\}$ to meromorphicity of $B$ there. 

In Section~\ref{s2}, we use the Geronimo-Case equations to prove Theorem~\ref{T1.5}. 
In Section~\ref{s3}, we use the second Szeg\H{o} map from OPRL to OPUC to prove 
Theorem~\ref{T1.6}. In Section~\ref{s4}, we extend the analysis of \cite{Jost2} to 
obtain Theorem~\ref{T1.7} from Theorem~\ref{T1.6}.  

\smallskip
This research was completed during my stay as a Lady Davis Visiting Professor 
at Hebrew University, Jerusalem. I'd like to thank H.~Farkas and Y.~Last for 
the hospitality of the Mathematics Institute at Hebrew University.

\section{The Geronimo-Case Equations and the $R^{-2}$ Result} \lb{s2} 

In this section, we will prove Theorem~\ref{T1.5} using a strategy similar to 
that used in \cite{MSF} to prove Theorem~\ref{T1.3}. There the critical element 
was the use of Szeg\H{o} recursion (1.13) and its $^*$, that is, 
\begin{equation} \lb{2.1} 
\Phi_{n+1}^* (z) =\Phi_n^*(z) -\alpha_n z \Phi_n(z) 
\end{equation} 
at $z$ and $1/\bar z$. 

Here we will instead use the Geronimo-Case equations \cite{GC} in the form introduced 
in \cite{Jost2}. Define 
\begin{equation} \lb{2.2} 
C_n (z) =z^n P_n \biggl( z + \f{1}{z}\biggr) 
\end{equation} 
The equations 
\begin{align} 
C_n(z) &= (z^2 -b_n z) C_{n-1}(z) + G_{n-1}(z) \lb{2.3}  \\
G_n(z) &= G_{n-1} (z) + [(1-a_n^2)z^2 - b_n z] C_{n-1}(z) \lb{2.4} 
\end{align} 
are the unnormalized GC equations. With initial condition $G_0(z)=C_0(z)=1$, they define 
monic polynomials of degree at most $n$. $C_n$ has the form \eqref{2.2}, and if  
\begin{equation} \lb{2.5} 
\sum_{n=1}^\infty\, \abs{a_n^2 -1} + \abs{b_n} <\infty 
\end{equation} 
then for $\abs{z}<1$, 
\begin{equation} \lb{2.6} 
\lim_{n\to\infty}\, G_n(z) = \biggl(\, \prod_{j=1}^\infty a_j\biggr) u(z) 
\end{equation} 
(see Theorem~A.3 of \cite{Jost2}). 
 
\eqref{2.3}/\eqref{2.4} have a structure somewhat like \eqref{1.13}/\eqref{2.1}. The 
difference is that \eqref{1.14} is replaced by 
\begin{equation} \lb{2.7} 
C_n(z) =z^{2n} C_n \biggl( \f{1}{z}\biggr) 
\end{equation} 
as is obvious from \eqref{2.2}. We introduce $f=\wti O(g)$ if and only if for all $\veps 
>0$, $\abs{f}/\abs{g}^{1-\veps} \to 0$. 

\begin{lemma} \lb{L2.1} If \eqref{1.21} holds, then for $z\in\bbD$,  
\begin{alignat}{3} 
&\text{\rm{(i)}} &\qquad \abs{G_n(z) -u(z)} &\leq \wti O(R^{-2n}) \lb{2.8} \\
&\text{\rm{(ii)}}&\qquad  \biggl|C_n(z) - \f{u(z)}{1-z^2}\biggr| 
&\leq \wti O([\max (\abs{z},R^{-1})]^{2n}) \lb{2.9} 
\end{alignat}  
\end{lemma} 

\begin{proof} (i) By Theorem~A.3 of \cite{Jost2}, 
\begin{equation} \lb{2.10} 
\lim_{n\to\infty}\, C_n(z) = \f{u(z)}{1-z^2} 
\end{equation}  
By \eqref{2.4} and $\sup_n \abs{C_n(z)} <\infty$, we see 
\begin{align*} 
\abs{G_n(z) -u(z)} &\leq \sum_{m=n}^\infty\, \abs{G_{m+1}(z) -G_m(z)} \\ 
&\leq \bigl(\, \sup_n \, \abs{C_n(z)}\bigr) \sum_{m=1}^\infty \, 
(\abs{ 1-a_{n+m}^2} + \abs{b_{n+m}}) \\ 
&= \wti O(R^{-2n}) 
\end{align*} 
since the series of bounds converges exponentially. 

\smallskip 
(ii) By \eqref{2.3}, 
\[
\abs{C_n -G_{n-1} -z^2 C_{n-1}} \leq \sup_n \, \abs{C_n(z)} \, \abs{b_n} 
\]
so iterating, 
\begin{align*} 
\biggl|C_n - \sum_{j=0}^{n-1} G_{n-j-1} z^{2j} \biggr| 
&\leq \abs{z}^{2n} + \sup_n\, \abs{C_n(z)} \sum_{j=0}^{n-1} \, \abs{b_{n-j}}\abs{z^{2j}} \\ 
&\leq \wti O(\max(\abs{z}, R^{-1})^{2n}) 
\end{align*} 
By \eqref{2.8}, 
\[
\biggl|\, \sum_{j=0}^{n-1}\, (G_{n-j-1}-u)z^{2j}\biggr| \leq 
\wti O(\max (\abs{z},R^{-1})^{2n}) 
\]
Since $\sum_j z^{2j} u=(1-z^2)^{-1}u$, we have \eqref{2.9}. 
\end{proof} 

\begin{proof}[Proof of Theorem~\ref{T1.5}] By \eqref{2.4} and \eqref{2.7} for 
$\abs{z} >1$, 
\begin{equation} \lb{2.11} 
\abs{G_{n+1} - G_n} \leq \biggl[\sup_n \, \biggl| C_n \biggl( \f{1}{z}\biggr)\biggr|\biggr] 
\abs{z}^{2n+2} [\abs{1-a_n^2} + \abs{b_n}]
\end{equation} 
which proves that for $1<\abs{z}<R$, $G_n$ converges uniformly, so by the maximum 
principle, we have convergence for $\abs{z}<R$, so $u$ has an analytic continuation 
to that region. In that region, 
\begin{align} 
u(z) &= 1+ \sum_{n=0}^\infty \, (G_{n+1}(z) -G_n(z)) \notag \\
&= 1+ \sum_{n=0}^\infty\, ((1-a_{n+1}^2) z^2 -b_{n+1}z) C_n(z) \lb{2.12} \\
&= \biggl( \f{u(\f{1}{z})}{1-\f{1}{z^2}}\biggr) (B(z)-1) + 1+ \sum_{n=0}^\infty 
f_n(z) \lb{2.13}   
\end{align} 
where 
\begin{equation} \lb{2.14} 
f_n(z) = ((1-a_{n+1}^2) z^2 -b_{n+1}z) z^{2n} \biggl( C_n \biggl(\f{1}{z}\biggr) 
- \f{u(\f{1}{z})}{1 - \f{1}{z^2}}\biggr)
\end{equation} 
Thus 
\begin{equation} \lb{2.15} 
(1-z^2) u(z) + u\biggl( \f{1}{z}\biggr) z^2 B(z) = 
u\biggl( \f{1}{z}\biggr) z^2 + (1-z^2) + \sum_{n=0}^\infty \, (1-z^2) f_n(z)
\end{equation} 

Each function $f_n$ is analytic in $\{z\mid\abs{z}>1\}$, so if we can prove that 
the sum converges uniformly in $\{z\mid 1<\abs{z} < R^2\}$, we know the left side 
of \eqref{2.15} has an analytic continuation in that region. 

By \eqref{2.9}, for $\abs{z}>1$, 
\[ 
\biggl| C_n\biggl( \f{1}{z}\biggr) -\f{u (\f{1}{z})}{1-\f{1}{z^2}}\biggr|
\leq \wti O\biggl( \max\biggl( \f{1}{\abs{z}}\, ,  R^{-1}\biggr)^{2n}\biggr) 
\]
so 
\[
\biggl| z^{2n} \biggl[ C_n \biggl( \f{1}{z}\biggr) - \f{u(\f{1}{z})}{1-\f{1}{z^2}}\biggr]\biggr| 
\leq \wti O(\max (1,\abs{z}R^{-1})^{2n}) 
\]
and thus, 
\[
\abs{f_n(z)} \leq \wti O(R^{-2n}) \wti O(\max (1,\abs{z} R^{-1})^{2n}) 
\]

For $1<\abs{z}<R$, this is $\wti O (R^{-2n})$ and so summable. For $R\leq\abs{z} <R^{-2}$, it is 
$\wti O((\abs{z}R^{-2})^{2n})$ and so also summable. 
\end{proof} 

If $u(\pm R^{-1})\neq 0$, \eqref{1.22} tells us that since $B$ has a singularity on the 
circle of radius $R$, so must $u$. However, if $u(R^{-1})=0$ and/or $u(-R^{-1})=0$, that 
zero can compensate for a pole in $B$ and $u$ can have a larger region of analyticity 
than $B$. This is exactly what happens in the case of noncanonical weights, as explained 
in \cite{Jost2}.

\section{The Second Szeg\H{o} Map and Jost Functions With No Bound States} \lb{s3} 

In \cite{Sz22a,Szb}, Szeg\H{o} defined two maps from the probability measures on 
$\partial\bbD$ invariant under $z\to\bar z$ to the probability measures on $[-2,2]$; let 
us call them $\Sz_1$ and $\Sz_2$. Both are injective, but only $\Sz_1$ is surjective  
--- and for this reason, $\Sz_1$ is the one most often used and studied (see 
\cite[Section~13.1]{OPUC2}). Here we will see that $\Sz_2$ is also exceedingly useful, 
especially for studying Jost functions analytic in a neighborhood of $\ol{\bbD}$ and 
nonvanishing on $\ol{\bbD}$ (i.e., $J$ has no bound states and no resonance at $\pm 2$). 

For a.c.\ measures, the relations are 
\begin{equation} \lb{3.1} 
d\mu = w(\theta)\, \f{d\theta}{2\pi} \qquad 
\Sz_1 (d\mu) = f_1 (x)\, dx \qquad 
\Sz_2 (d\mu) = f_2(x)\, dx 
\end{equation} 
where $w(\theta) =w(-\theta)$ and (formulae (13.1.6) and (13.2.22) of \cite{OPUC2}) 
\begin{align} 
f_1(x) &= \pi^{-1} (4-x^2)^{-1/2} w\biggl(\arccos \biggl( \f{x}{2}\biggr)\biggr) \lb{3.2}\\ 
f_2(x) &= \pi^{-1} c^2 (4-x^2)^{1/2} w\biggl( \arccos\biggl(\f{x}{2}\biggr)\biggr) \lb{3.3}  
\end{align} 
where 
\begin{equation} \lb{3.4} 
c=[2(1-\abs{\alpha_0}^2)(1-\alpha_1)]^{-1/2} 
\end{equation} 
Taking into account that $\Sz_1$ is a bijection of even measures on $\partial\bbD$ and all 
measures on $[-2,2]$, we see that 
\begin{equation} \lb{3.5} 
d\gamma\in\ran(\Sz_2)\Leftrightarrow \int_{-2}^2 (4-x^2)^{-1}\, d\gamma(x) < \infty 
\end{equation} 

\begin{proposition}\lb{P3.1} If $d\gamma$ has a Jost function $u$ analytic in a neighborhood 
of $\ol{\bbD}$ and nonvanishing on $\ol{\bbD}$, then $d\gamma \in\ran(\Sz_2)$. 
\end{proposition} 

\begin{proof} Since $d(2\cos\theta)=-2\sin\theta\,d\theta$ and $(4-2\cos^2\theta)= 
4\sin^2\theta$, by \eqref{1.4} the right side of \eqref{3.5} is equivalent to 
\begin{align*} 
\f{1}{\pi} \int_0^\pi \f{d\theta}{\abs{u(e^{i\theta})}^2} &= 
\int_0^\pi (\sin^2 \theta)^{-1} f(2\cos\theta)\sin\theta\, d\theta \\
&= 2 \int_{-2}^2 (4-x^2)^{-1} f(x)\,dx <\infty 
\end{align*} 
which is true if $\abs{u}$ is bounded away from zero. 
\end{proof} 

For our purposes, what is critical is: 

\begin{theorem}\lb{T3.2} Let $d\mu$ be a nontrivial probability measure on $\partial\bbD$ 
obeying the Szeg\H{o} condition with Verblunsky coefficients $\{\alpha_n\}_{n=0}^\infty$ 
and Szeg\H{o} function $D(z_0)$. Let $d\gamma=\Sz_2 (d\mu)$ have Jost function $u$ and 
Jacobi parameters $\{a_n,b_n\}_{n=1}^\infty$. Then 
\begin{align} 
b_{n+1} &= \alpha_{2n} -\alpha_{2n+2} - \alpha_{2n+1} (\alpha_{2n} + \alpha_{2n+2}) \lb{3.6} \\
a_{n+1}^2-1 &= \alpha_{2n+1} -\alpha_{2n+3} - \alpha_{2n+2}^2 (1-\alpha_{2n+3})(1+\alpha_{2n+1}) 
- \alpha_{2n+3}\alpha_{2n+1} \lb{3.7} \\
u(z) &= (1-\abs{\alpha_0}^2)(1-\alpha_1) D(z)^{-1} \lb{3.8}
\end{align} 
\end{theorem} 

\begin{remark} Formulae of the form \eqref{3.6}/\eqref{3.7} for $\Sz_1$ go back to Geronimus 
\cite{Ger46,GBk1}. For $\Sz_2$, the earliest reference I am aware of is Berriochoa, Cachafeiro, 
and Garc\'ia-Amor \cite{BCG1}; see also \cite{KilNen}. 
\end{remark} 

\begin{proof} \eqref{3.6}/\eqref{3.7} are (13.2.20)/(13.2.21) of \cite{OPUC2}. To see 
\eqref{3.8}, note that, by \eqref{1.4} and \eqref{3.3}, 
\begin{align*} 
\abs{u(e^{i\theta})}^{-2} &= \pi f_2 (2\cos\theta)(\sin\theta)^{-1} \\
&= 2c^2 w \\
&= 2c^2 \abs{D}^2 
\end{align*} 
Thus, the absolute value of \eqref{3.8} holds if $\abs{z}=e^{i\theta}$. Since both 
sides are analytic, nonvanishing on  $\bbD$, and positive at $z=0$, \eqref{3.8} holds 
for all $z$. 
\end{proof} 

\eqref{3.6}/\eqref{3.7} first of all provide a second proof of Theorem~\ref{T1.5} 
in case $u$ is nonvanishing on $\ol{\bbD}$ and, more importantly, show generically 
that $R^2$ is optimal. We note first: 

\begin{proposition}\lb{P3.3} We have that 
\begin{equation} \lb{3.9} 
B(z) =\alpha_0 z^{-1} + \alpha_1 +1 + (S(z)-1)(1-z^{-2})+Q(z)
\end{equation}  
where, if 
\begin{equation} \lb{3.10} 
\limsup \abs{\alpha_n}^{1/n}=R^{-1} 
\end{equation} 
then $Q$ is analytic in $\{z\mid\abs{z}<R^2\}$. 
\end{proposition} 

\begin{proof} By \eqref{1.9}, \eqref{1.15}, and \eqref{3.6}/\eqref{3.7}, we have \eqref{3.9} 
where 
\begin{equation} \lb{3.11} 
\begin{split}
Q(z) =\sum_{n=0}^\infty &\alpha_{2n+1} (\alpha_{2n}+\alpha_{2n+2}) z^{2n+1} \\ 
& +  \{\alpha_{2n+2}^2 (1-\alpha_{2n+3})(1+\alpha_{2n+1}) + \alpha_{2n+3}\alpha_{2n+1}\} 
z^{2n+2} 
\end{split}
\end{equation} 
By \eqref{3.10}, $Q(z)$ is analytic in $\{z\mid \abs{z} <R^2\}$. 
\end{proof} 

\begin{proof}[Second proof of Theorem~\ref{T1.5} when $u$ is nonvanishing 
on $\ol{\bbD}$] As we will show below (see Lemma~\ref{L3.5}), \eqref{1.26} 
implies \eqref{3.10}. By Theorem~\ref{T1.3} and \eqref{3.8}, we conclude that 
\begin{equation} \lb{3.12} 
(z^2 -1) \bigl[ u(z) - \ol{u(1/\bar z)}\, S(z)\bigr] 
\end{equation} 
is analytic in $\{z\mid R^{-1} < \abs{z} < R^3\}$. By Proposition~\ref{P3.3}, 
\begin{equation} \lb{3.13} 
z^2 B(z) - (z^2-1)S(z) 
\end{equation}
is analytic in $\{z\mid \abs{z}<R^2\}$, so by \eqref{3.12}, the function in \eqref{1.22} 
is analytic in $\{z\mid R^{-1} < \abs{z} <R^2\}$. 
\end{proof} 

\begin{example}\lb{E3.4} Suppose $\alpha_{2n}\equiv 0$ (true if and only if $b_n 
\equiv 0$) and $\alpha_{2n+1}=R^{-(2n+1)}$. Then, by \eqref{3.10}, 
\begin{align*} 
Q(z) &=\sum_{n=0}^\infty z^{2n+2} R^{-4n-4} \\
&= z^2 R^{-4} (1-z^2 R^{-4})^{-1} 
\end{align*} 
has poles at $z=\pm R^2$. This shows that \eqref{1.22} may not be analytic in any larger 
annulus than $\{z\mid R^{-1} <\abs{z} <R^2\}$. It is also clear that by a similar analysis,  
if $B$ is meromorphic in $\{z\mid\abs{z}<R^{1+\veps}\}$, then generically \eqref{1.22} will 
have singularities on the circle of radius $R^2$. The change from $R^3$ to $R^2$ in going 
from Theorem~\ref{T1.3} to Theorem~\ref{T1.5} is due to the quadratic terms in 
\eqref{3.6}/\eqref{3.7}. 
\qed
\end{example} 

Above we used and below we will need: 

\begin{lemma}\lb{L3.5} Suppose $\{\alpha_n\}_{n=0}^\infty$ and $\{a_n,b_n\}_{n=1}^\infty$ are 
related by \eqref{3.6}/\eqref{3.7}, and  
\begin{equation} \lb{3.14} 
\limsup \abs{\alpha_n}^{1/n} \equiv R_1^{-1} <1 \qquad 
\limsup_{n\to\infty}\, (\abs{a_n^2 -1} +\abs{b_n})^{1/2n} \equiv R_2^{-1} <1
\end{equation} 
Then $R_1=R_2$. Moreover, $\alpha_n$ has a complete asymptotic series if and only if 
$\{a_n,b_n\}$ do, and if $T$ is the set of powers that enter for $\{a_n,b_n\}_{n=1}^\infty$ 
{\rm{(}}i.e., $T$ are the poles of $B${\rm{)}} and $\wti T$ for $\{\alpha_n\}_{n=0}^\infty$ 
{\rm{(}}i.e., $\wti T$ are the poles of $S${\rm{)}}, then 
\begin{equation} \lb{3.15} 
T\subset\wti\bbG(\wti T) \qquad \wti T\subset\wti\bbG(T) 
\end{equation} 
\end{lemma} 

\begin{remark} For $\Sz_1$, there are equations similar to \eqref{3.6}/\eqref{3.7} which 
have solutions where $\{a_n,b_n\}_{n=0}^\infty$ has rapid decay while $\alpha_{2n+1}\sim 
n^{-1}$ at infinity. (Indeed, for $\Sz_1$ but {\it not\/} $\Sz_2$, this happens for $J_0$; 
see Example~13.1.3 revisited in \cite{OPUC2}.) In fact, the results in this paper plus 
\cite{MSF} imply that $R_1^{-1}<1$ if and only if $R_2^{-1} <1$. 
\end{remark} 

\begin{proof} It follows from \eqref{3.6} that if $R_1^{-1}, R_2^{-1} <1$, then 
$R_2 =R_1$ since the nonleading terms are exponentially small. In addition, 
if $\alpha_n$ has a complete asymptotic series, one gets that $b_{n+1}$ and 
$a_{n+1}^2 -1$ individually have asymptotic series in $\mu_k^{-2n}$ with 
$\mu_k\in \cup_{j=1}^\infty \bbG^{(j)}(\wti T)$. Since 
\begin{align} 
c_1\mu_k^{-2n} &=\tfrac12\, (c_1+c_2) \mu_k^{-2n} + 
\tfrac12\, (c_1-c_2)(-\mu_k)^{-2n} \lb{3.16} \\ 
c_2\mu_k^{-2n-1} &=\tfrac12\, (c_1+c_2) \mu_k^{-2n-1} + 
\tfrac12\, (c_1-c_2)(-\mu_k)^{-2n-1} \lb{3.17} 
\end{align} 
we can combine into a single series by taking $-\mu$'s as well as $\mu$'s. 

For the converse, note that since the $\alpha$'s decay exponentially, 
\[
\alpha_{2n} =\sum_{m=0}^\infty b_{n+m+1} + O(R^{-2n}) 
\]
and similarly for $\alpha_{2n+1}$ and $\sum_{m=0}^\infty (a_{n+m+1}^2 -1)$. Plugging this 
into \eqref{3.6} and summing yields $\alpha_{2n}$ and $\alpha_{2n+1}$ as explicit sums of 
products of four or fewer $b$'s and $(1-a^2)$'s plus an error of $O(R^{-3n})$. Iterating 
gives explicit formulae for $\alpha$'s as ``polynomials" in $b$ and $1-a^2$ of degree 
$k$ plus an error of order $O(R^{-(k+2)n})$. This shows that if $a$ and $b$ have 
asymptotic series to order $R^{-(k+1)n}$, so do $\alpha_{2n}$ and $\alpha_{2n-1}$ 
with rates in $\cup_{j=1}^\infty \bbG^{(j)}(T)$. Using formulae like \eqref{3.16}/\eqref{3.17}, 
we can combine to a single series by using $-\mu$'s, so $\wti T\subset\wti\bbG(T)$. 
\end{proof} 

\begin{proof}[Proof of Theorem~\ref{T1.6} and Theorem~\ref{T1.2} when $u$ is nonvanishing 
on $\ol{\bbD}$] Since $u$ is nonvanishing on $\ol{\bbD}$, $\gamma\in\ran(\Sz_2)$, so 
we can define $S,\alpha_n$, etc. If $u$ is entire meromorphic, by \eqref{3.8}, so is $D^{-1}$. 
Thus, by Theorem~\ref{T1.4}, $S$ is entire meromorphic, and if $\wti T$ are the poles of 
$S$, then 
\[
\wti T\subset\bbG(P) 
\]
By Lemma~\ref{L3.5}, $B$ is meromorphic and 
\[
T\subset\wti\bbG(\wti T)\subset\wti\bbG(\bbG(P)) =\wti\bbG(P) 
\]

Conversely, if $B$ is entire meromorphic, by Lemma~\ref{L3.5}, so is $S$, and if 
$\wti T$ are the poles of $S$, then 
\[
\wti T \subset\wti\bbG(T) 
\]
By Theorem~\ref{T1.4}, $D^{-1}$, and so $u$, is entire meromorphic and  
\[
P\subset\bbG(\wti T)\subset\bbG(\wti\bbG(T)) = \wti\bbG(T) 
\qedhere 
\] 
\end{proof}

\section{Coefficient Stripping and Jost Functions With Bound States} \lb{s4} 

As in \cite{Jost2}, we will go from the no bound state theorem to the general case 
(i.e., in our situation, Theorem~\ref{T1.6} to Theorem~\ref{T1.7}) by coefficient 
stripping, that is, pass from $J$ to the Jacobi matrix $J^{(m)}$ with Jacobi parameters 
$\{a_{n+m}, b_{n+m}\}_{n=1}^\infty$. By Theorem~3.1 of \cite{Jost1}, if $J$ has 
a Jost function analytic in a neighborhood of $\ol{\bbD}$, there exists a $k$ with 
$\sigma (J^{(k)})=[-2,2]$, and by a slight extension of the argument, we can also 
suppose its Jacobi function obeys $u^{(k)}(\pm 1)\neq 0$ (for if $\sigma(J^{(k-1)})= 
[-2,2]$, if $u^{(k+1)}$ vanishes at $\pm 1$, $M^{(k-1)}(z)$ has a pole there and 
$u^{(k)}=u^{(k-1)} M^{(k-1)}$ is nonvanishing). Thus, we claim that we need only prove 
(as we shall do below) that 

\begin{theorem}\lb{T4.1} If $P=P_1\cup P_2$ as in Theorem~\ref{T1.7} and we make the 
$J$-dependence explicit, then 
\begin{equation} \lb{4.1} 
P(J)=P(J^{(1)}) 
\end{equation} 
\end{theorem} 

\begin{proof}[Proof of Theorems~\ref{T1.1} and \ref{T1.7} given Theorems~\ref{T4.1} 
and \ref{T1.6}] Theorem~\ref{T1.7} implies Theorem~\ref{T1.1}. By \eqref{4.1} and 
induction, $P(J)=P(J^{(k)})$ where $k$ is chosen as above so Theorem~\ref{T1.6} is 
applicable. \eqref{1.26} for $J^{(k)}$ implies it for $J$\!. 
\end{proof} 

As in \cite{Jost2}, we will make use of the $M$-function, its connection to $u$, and 
the update relations. We define (consistently with \eqref{1.5})
\begin{equation} \lb{4.2} 
M^{(k)}(z)=\langle\delta_1, (z+z^{-1} -J^{(k)})^{-1}\delta_1\rangle 
\end{equation} 
for $z\in\bbD\backslash \{w\mid w+w^{-1} \in\sigma(J^{(k)})\}$. $M^{(k)}$ has poles 
at the set in $\bbD$ with $w+w^{-1}\in\sigma (J^{(k)})$, and $u^{(k)}$ has zeros there. 
The update equations ((2.4)/(2.5) of \cite{Jost2}) are (initially for $z\in\bbD$) 
\begin{align} 
u^{(k+1)}(z) &=a_{k+1} z^{-1} u^{(k)}(z) M^{(k)}(z) \lb{4.3} \\
M^{(k)}(z)^{-1} &= z+z^{-1} -b_{k+1} -a_{k+1}^2 M^{(k+1)}(z) \lb{4.4}  
\end{align} 
Moreover, we have the analytic continuation of \eqref{1.4} plus $\pi f(2\cos\theta) 
=\Ima M(e^{i\theta})$ for $\theta\in [0,\pi]$, 
\begin{equation} \lb{4.5} 
[M(z)-\ol{M}(1/\bar z)]\, \ol{u(1/\bar z)}\, u(z) = z-z^{-1} 
\end{equation} 

\eqref{4.5} can be used to meromorphically continue $M$ from $\bbD$ to $\bbC$ if $u$ is 
entire meromorphic. Once one makes these continuations, \eqref{4.3}/\eqref{4.4} extend 
to all $z\in\bbC$ (as meromorphic relations including possible cancellations of poles 
and zeros). \eqref{4.3}/\eqref{4.4} also show that if $u$ is entire meromorphic, so 
is $u^{(1)}$. 

We begin by rephrasing the set $P_2$: 

\begin{proposition}\lb{P4.2} $z_0\in\bbC\backslash\ol{\bbD}$ is in $P_2$ if and only 
if 
\begin{SL} 
\item[{\rm{(i)}}] $z_0$ is not a pole of $u$. 
\item[{\rm{(ii)}}] Both $z_0$ and $z_0^{-1}$ are poles of $M$\!. 
\end{SL} 
\end{proposition} 

\begin{remark} In $\bbD$, all poles of $M$ are real, so (ii) implies that $z_0$ 
is real. 
\end{remark} 

\begin{proof} By definition, $z_0\in P_2$ if and only if $u(z_0^{-1})=0$, $z_0$ is 
not a pole of $u$, and \eqref{1.27} holds. Since $z_0^{-1}\in\bbD$, $u(z_0^{-1})=0 
\Leftrightarrow z_0 + z_0^{-1}\in\sigma(J)\Leftrightarrow z_0^{-1}$ is a pole of 
$M(z)$. 

As shown in \cite{Jost2}, by \eqref{4.5}, if $u(z_0)=0$, $z_0$ has a pole of  
$M$ of order two or more and, of course, \eqref{1.27} holds at $z_0^{-1}$ since 
the left side is infinite. If $u(z_0)\neq 0$, \eqref{1.27} is precisely the condition, 
via \eqref{4.5}, that $M(z)$ has a pole at $z_0$. 
\end{proof} 

We have been careful in considering situations where $z_0$ is a pole of $u$ 
and $z_0^{-1}$ is a zero of $u$. We need to consider that case separately: 

\begin{proposition}\lb{P4.3} If $z_0\in\bbC\backslash\bbD$ is a pole of $u$ and 
$z_0^{-1}$ is a zero of $u$, then $z_0$ is a pole of $u^{(1)}$. 
\end{proposition} 

\begin{proof} Consider \eqref{4.5} near $z=z_0$. Zeros of $u$ in $\bbD$ are 
simple, so $u(z)\, \ol{u(1/\bar z)}$ either has a pole at $z_0$ or a finite 
nonzero limit. Thus, \eqref{4.5} shows $M(z)-\ol{M(1/\bar z)}$ must be regular 
(perhaps even zero) at $z_0$. Since $\ol{M(1/\bar z)}$ has a pole at $z_0$, 
$M(z)$ must have a pole there also. It follows that $u^{(1)}=uM$ has a pole    
(indeed, at least a second-order pole) at $z_0$. 
\end{proof} 

\begin{proposition}\lb{P4.4} If $z_0\in P_2 (J)$, then $z_0\in P_1(J^{(1)})$. 
\end{proposition} 

\begin{remarks} 1. For $z_0$ within the disk of analyticity of $u$, this result 
is in \cite{Jost2}. The proof here is essentially identical. 

2. $z_0\in P_2(J)$ is essentially a statement of the vanishing of a ``resonance 
eigenfunction," so this says that such eigenfunctions cannot have successive 
zeros because of a second-order difference equation. 
\end{remarks} 

\begin{proof} Suppose first that $u(z_0)\neq 0$. By Proposition~\ref{P4.2}, $M$ has 
a pole at $z_0$, so $u^{(1)}=uM$ has a pole at $z_0$. 

If $u$ has a $k$th-order zero, $k\geq 1$, at $z_0$, $\ol{u(1/\bar z)}\,u(z)$ has a 
$(k+1)$st-order zero, so $M(z)-\ol{M(1/\bar z)}$ has a $(k+1)$st-order pole $z_0$ 
by \eqref{4.5}. Since $M$ has simple poles at points in $\bbD$ like $1/z_0$, $M$ 
has to have a $(k+1)$st-order pole at $z_0$. Thus, $u^{(1)}=uM$ has a pole at 
$z_0$. 
\end{proof} 

\begin{proposition}\lb{P4.5} If $z_0\in P_1(J)$ and $z_0\notin P_1(J^{(1)})$, then 
$z_0\in P_2(J^{(1)})$. 
\end{proposition} 

\begin{proof} By \eqref{4.4}, poles of $M^{(1)}(z)$ are precisely at zeros of $M(z)$. 
Thus, by Proposition~\ref{P4.2}, we need to prove that  
\begin{equation} \lb{4.6} 
z_0\in P_1(J) + z_0\notin P_1(J^{(1)}) \Rightarrow M(z_0) = 
M\biggl(\f{1}{z_0}\biggr) =0 
\end{equation} 

Since $u$ has a pole at $z_0$ and $u^{(1)}=uM$ does not, $M(z_0)=0$. By 
Proposition~\ref{P4.3}, $z_0\notin P_1 (J^{(1)})$ implies $u(1/\bar z) \neq 0$. 
Thus, $u(z)\, \ol{u(1/\bar z)}$ has a pole at $z_0$. \eqref{4.5} then implies 
that $M(z) - \ol{M(1/\bar z)}|_{z=z_0}=0$. Since $M(z_0) =0$, 
we conclude $M(1/\bar z_0)=0$. This proves \eqref{4.6}. 
\end{proof} 

We also need some results that go back from $z_0\in P(J^{(1)})$. 

\begin{proposition}\lb{P4.6} If $z_0\in P_2 (J^{(1)})$, then $z_0\in P_1 (J)$. 
\end{proposition} 

\begin{proof} By Proposition~\ref{P4.2}, $z_0$ and $z_0^{-1}$ are poles of $M^{(1)}$, 
so by \eqref{4.4}, they are zeros of $M$\!. As in the proof of Proposition~\ref{P4.4}, 
if $u^{(1)}$ has a $k$th-order zero (including $k=0$, i.e., $u^{(1)}(z_0)\neq 0$), 
then $M^{(1)}(z)$ has a $(k+1)$st-order pole there, and so $M(z)$ has a $(k+1)$st-order 
zero. This is only consistent with $u^{(1)}=uM$ if $u$ has a pole at $z_0$.  
\end{proof} 

\begin{proposition}\lb{P4.7} If $z_0\in P_1 (J^{(1)})$ and $z_0\notin P_1(J)$, 
then $z_0\in P_2(J)$. 
\end{proposition} 

\begin{proof} By hypothesis, $z_0$ is not a pole of $u$, and so (i) of 
Proposition~\ref{P4.2} holds. So we need only show that $M(z)$ has poles at $z_0$ 
and $z_0^{-1}$. Suppose $u$ has a $k$th-order zero at $z_0$ (including $k=0$, i.e., 
$u(z_0)\neq 0$). By \eqref{4.3} and the fact that $z_0$ is a pole of $u^{(1)}$, 
we see that $z_0$ is a $(k+1)$st-order pole of $M(z)$ and, in particular, a pole of 
$M(z)$ (since $k+1\geq 1$). 

If $k\geq 1$, this is only consistent with \eqref{4.5} if $\ol{u(1/\bar z)}$ has 
a zero at $z_0$ since the possible pole of $M$ at $z_0^{-1}$ is of order $1$ and 
cannot cancel the $(k+1)$st-order pole at $z_0$. Thus, $\ol{u(1/\bar z_0)} =0$, so 
$z_0^{-1}$ is real and a pole of $M(z)$, that is, (ii) of Proposition~\ref{P4.2} 
holds and $z_0\in P_2(J)$. 

If $k=1$ and $M(z)$ does not have a pole at $\bar z_0^{-1}$, then $M(z) - \ol{M} 
(1/\bar z)$ has a pole at $z_0$, while $\ol{u(1/\bar z_0)}\neq 0\neq u(z_0)$ 
(since $k=1$ and $1/\bar z_0$ is not a pole), violating \eqref{4.5}. Thus, $M$ 
must have a pole at $z_0$ and $z_0\in P_2(J)$ by Proposition~\ref{P4.2}. 
\end{proof} 

\begin{proof}[Proof of Theorem~\ref{T4.1}] If $z_0\in P(J^{(1)})$, either $z_0\in 
P_2(J^{(1)})\Rightarrow z_0\in P_1(J)$ (by Proposition~\ref{P4.6}) or $z_0\in 
P_1(J^{(1)})\Rightarrow z_0\in P_1(J)\cup P_2(J)$ (by Proposition~\ref{P4.7}). 
Thus, $P(J^{(1)})\subset P(J)$. 

If $z_0\in P(J)$, either $z_0\in P_2(J)\Rightarrow z_0\in P_1 (J^{(1)})$ 
(by Proposition~\ref{P4.4}) or $z_0\in P_1(J)\Rightarrow z_0\in P_1 (J^{(1)}) 
\cup P_2 (J^{(1)})$ (by Proposition~\ref{P4.5}). Thus, $P(J)\subset P(J^{(1)})$. 
\end{proof}

\bigskip

\end{document}